\documentclass[reqno]{amsart}
\usepackage{amssymb}
\usepackage{verbatim}

\begin{document}

\title[(Anti)symmetric multivariate exponential
functions]
{(Anti)symmetric multivariate exponential
functions and corresponding Fourier transforms}

\author{A.~Klimyk}
\address{Bogolyubov Institute for Theoretical Physics,
         Kiev 03143, Ukraine}
 \email{aklimyk@bitp.kiev.ua}

\author{J.~Patera}
\address{Centre de Recherches Math\'ematiques,
         Universit\'e de Montr\'eal,
         C.P.6128-Centre ville,
         Montr\'eal, H3C\,3J7, Qu\'ebec, Canada}
\email{patera@crm.umontreal.ca}

 \begin{abstract}
We define and study symmetrized and antisymmetrized
multivariate exponential functions. They are defined as
determinants and antideterminants of matrices whose
entries are exponential functions of one variable.
These functions are eigenfunctions of the Laplace
operator on corresponding fundamental domains
satisfying certain boundary conditions.
To symmetric and antisymmetric multivariate exponential functions
there correspond Fourier transforms. There are three
types of such Fourier transforms: expansions into
corresponding Fourier series, integral Fourier transforms, and
multivariate finite Fourier transforms.
Eigenfunctions of the integral Fourier transforms are found.

 \end{abstract}

\maketitle

\section{Introduction}
In mathematical and theoretical physics, very often we deal with functions
on the Euclidean space $E_n$ which are symmetric or antisymmetric with
respect to the permutation (symmetric) group $S_n$. For example, such
functions describe collections of identical particles.
Symmetric and antisymmetric solutions appear in the
theory of integrable systems. Characters of finite dimensional
representations of semisimple Lie algebras are symmetric functions.
Moreover, according to the Weyl formula for these characters,
each such character is a ratio of antisymmetric functions.

The aim of this paper is to describe and to study
symmetrized and antisymmetrized multivariate exponential functions
and the corresponding Fourier transforms.
Antisymmetric multivariate exponential functions (we denote them by
$E^-_\lambda(x)$, $\lambda=(\lambda_1,\lambda_2,
\dots,\lambda_n)$, $x=(x_1,x_2,\dots,x_n)$) are determinants of
$n\times n$ matrices, whose entries are the usual exponential
functions of one variable, $E^-_\lambda(x)=
\det (e^{2\pi{\rm i}\lambda_ix_j})_{i,j=1}^n$.
Symmetric multivariate exponential functions
$E^+_\lambda(x)$
are antideterminants of the same $n\times n$ matrices (a definition of
antideterminants see below).

As in the case of the exponential functions of one variable, we may
consider three types of antisymmetric and symmetric multivariate exponential
functions:
 \medskip

(a) functions $E^-_m(x)$ and
$E^+_m(x)$ with $m=(m_1,m_2,
\dots,m_n)$, $m_i\in \mathbb{Z}$; they determine Fourier
series expansions in multivariate symmetric and antisymmetric exponential
functions;
 \medskip

(b) functions $E^-_\lambda(x)$ and
$E^+_\lambda(x)$ with $\lambda=(\lambda_1,\lambda_2,
\dots,\lambda_n)$, $\lambda_i\in \mathbb{R}$; these functions determine
integral multivariate Fourier transforms;
 \medskip

(c) functions $E^-_\lambda(x)$ and
$E^+_\lambda(x)$, where $x=(x_1,x_2,
\dots,x_n)$ take a finite set of values; they determine multivariate
finite Fourier transforms.
 \medskip

 Functions (b) are antisymmetric (symmetric) with respect to elements
of the permutation group $S_n$. (Anti)symmetries of functions (a)
are described by a wider group, since exponential functions
$e^{2\pi{\rm i}mx}$, $m\in \mathbb{Z}$, of one variable
are invariant with respect to shifts $x\to x+k$, $k\in \mathbb{Z}$.
These (anti)symmetries are described by elements of the
affine symmetric group $S_n^{\rm aff}$ which is a product of the
group $S_n$ and the group $T_n$, consisting of shifts in the space $E_n$ by
vectors $r=(r_1,r_2,\dots,r_n)$, $r_j\in \mathbb{Z}$.
A fundamental domain $F(S_n^{\rm aff})$ of the group $S_n^{\rm aff}$
is a certain bounded subset of $\mathbb{R}^n$ (see subsection 2.2).

The functions
$E^+_\lambda(x)$ give solutions of the Neumann
boundary value problem on a closure of the fundamental domain $F(S_n)$.
The functions
$E^-_\lambda(x)$ are solutions of
the Laplace equation $\Delta f=\mu f$ on the domain
$F(S_n)$ vanishing on the boundary $\partial F(S_n)$
of $F(S_n)$.

Functions on the fundamental domain $F(S_n^{\rm aff})$ can be expanded into
series in the functions (a). These expansions are an analogue of the
usual Fourier series for functions of one variable. Functions (b)
determine an (anti)symmetrized Fourier integral transforms on the
fundamental domain $F(S_n)$ of the symmetric group $S_n$. This
domain consists of points $x\in E_n$ such that $x_1>x_2>\cdots >x_n$.

Functions (c) are used to determine (anti)symmetric finite (that is, on a finite set)
Fourier transforms. These finite Fourier transforms are given on grids
consisting of points in the fundamental domain $F(S_n^{\rm aff})$.

Symmetric and antisymmetric exponential functions are closely related to
symmetric and antisymmetric orbit functions defined in
\cite{P04}, \cite{P-SIG-05} and studied in detail in \cite{KP06}
and \cite{KP07}. In fact, symmetric and antisymmetric exponential functions
are connected with orbit functions corresponding to the Coxeter--Dynkin
diagram $A_n$. Discrete orbit function transforms, corresponding to
Coxeter--Dynkin diagrams of low order, were detailly studied and it was
shown that they are very useful for applications \cite{AP}--\cite{Appl-5}.

The exposition of the theory of orbit functions in \cite{KP06}
and \cite{KP07} strongly depends on the theory of Weyl groups,
properties of root systems, etc. In this paper we avoid this dependence. We
use only the permutation (symmetric) group and properties of
determinants and antideterminants.
It is well-known that a determinant $\det (a_{ij})_{i,j=1}^n$ of the
$n\times n$ matrix $(a_{ij})_{i,j=1}^n$ is defined as
$$
\det (a_{ij})_{i,j=1}^n=\sum_{w\in S_n} (\det w)
a_{1,w(1)}a_{2,w(2)}\cdots a_{n,w(n)}
$$
\begin{equation}  \label{det}
= \sum_{w\in S_n} (\det w)
a_{w(1),1}a_{w(2),2}\cdots a_{w(n),n},
\end{equation}
where $S_n$ is the symmetric group of $n$ symbols
$1,2,\dots,n$, the set $(w(1),w(2),\dots$, $w(n))$ means the set
$w(1,2,\dots,n)$, and  $\det w$ denotes a determinant of the transform $w$, that is,
$\det w=1$ if $w$ is an even permutation and $\det w=-1$ otherwise.
Along with a determinant, we shall use an antideterminant $\det^+$ of
the matrix $(a_{ij})_{i,j=1}^n$ which is defined as a sum of all
terms, entering to the expression for the corresponding determinant,
taken with the sign +,
\begin{equation}  \label{dett}
{\det}^+ (a_{ij})_{i,j=1}^n=\sum_{w\in S_n} a_{1,w(1)}a_{2,w(2)}\cdots
a_{n,w(n)}= \sum_{w\in S_n}
a_{w(1),1}a_{w(2),2}\cdots a_{w(n),n}.
\end{equation}

Symmetrized and antisymmetrized multivariate polynomials were studied by
several authors (see, for example, \cite{KMc62} and \cite{BSX95}). In this
paper we investigate symmetric and antisymmetric multivariate exponential
functions.

\section{Symmetric and antisymmetric multivariate\\ exponential
functions}

\subsection{Definition}
A symmetric multivariate exponential function of $x=(x_1,x_2$, $\dots,x_n)$ is
defined as the function
\begin{alignat}{2}  \label{det-SS}
E_\lambda^+(x)\equiv &
E^+_{(\lambda_{1},\lambda_2,\dots ,\lambda_{n})}(x)=
  {\det}^+ \left( e^{2\pi{\rm
i}\lambda_ix_j}\right)_{i,j=1}^{n}
\notag\\
= & { {\det}^+ \left(
 \begin{array}{cccc}
 e^{2\pi{\rm i}\lambda_1x_1}& e^{2\pi{\rm i}\lambda_1x_2}&\cdots & e^{2\pi{\rm
i}\lambda_1x_{n}}\\
  e^{2\pi{\rm i}\lambda_2x_1}& e^{2\pi{\rm i}\lambda_2x_2}&\cdots & e^{2\pi{\rm
i}\lambda_2x_{n}}\\
  \cdots & \cdots &
\cdots & \cdots \\
 e^{2\pi{\rm i}\lambda_{n}x_1}& e^{2\pi{\rm i}\lambda_{n}x_2}&\cdots & e^{2\pi{\rm
i}\lambda_{n}x_{n}} \end{array} \right) }
\notag\\
\equiv & \sum_{w\in S_n}
 e^{2\pi{\rm i}\lambda_1x_{w(1)}} e^{2\pi{\rm i}\lambda_2x_{w(2)}}\cdots  e^{2\pi{\rm
i}\lambda_nx_{w(n)}}
= \sum_{w\in S_n} e^{2\pi{\rm i}\langle \lambda, wx \rangle} ,
 \end{alignat}
where $\lambda=(\lambda_1,\lambda_2,\dots,\lambda_n)$ is a set of real
numbers, which determines the function $E_\lambda^+(x)$,
and $\langle \lambda,x  \rangle$ denotes the scalar product in the
$n$-dimensional Euclidean space $E_n$,
$\langle \lambda,x  \rangle=\sum_{i=1}^n\lambda_ix_i$.
When $\lambda_1,\lambda_2,\dots,\lambda_n$ are
integers, we denote this set of numbers as $m\equiv (m_1,m_2,\dots,m_n)$,
\begin{equation}   \label{detI}
E_{m_1,m_2,\dots,m_n}^+(x)={\det}^+ \left( e^{2\pi{\rm
i}m_ix_j}\right)_{i,j=1}^{n}.
 \end{equation}

It is seen from the expression \eqref{dett} for an antideterminant
${\det}^+$ that the symmetric exponential functions
$E^+_\lambda(x)$ satisfy the relation
\begin{equation}  \label{x_1+x_2}
E^+_\lambda(x_1+a,x_2+a,\dots,x_n+a)=
e^{2\pi{\rm i}(\lambda_1+\lambda_2+\cdots +\lambda_n)a}
E^+_\lambda(x).
\end{equation}
This means that it is enough to consider the function
$E^+_\lambda(x)$ on the hyperplane
\[
x_1+x_2+\cdots +x_n=b,
\]
where $b$ is a fixed number (we denote this hyperplane
by $\mathcal{H}_b$). A transition from one hyperplane
$\mathcal{H}_b$ to another $\mathcal{H}_c$ is fulfilled by
multiplication by a usual exponential function
$e^{2\pi{\rm i}\vert\lambda \vert (c-b)}$, where
$\vert \lambda \vert=\lambda_1+\lambda_2+\cdots +\lambda_n$,
\[
E^+_\lambda(x)\vert_{x\in {\mathcal{H}_b}}=
e^{2\pi{\rm i}(\lambda_1+\lambda_2+\cdots +\lambda_n)(c-b)}
E^+_\lambda(x)\vert_{x\in {\mathcal{H}_c}}.
\]
It is useful to consider the functions
$E^+_\lambda(x)$ on the hyperplane
$\mathcal{H}_0$.
For  $x\in\mathcal{H}_0$ we have the relation
\begin{equation}  \label{x_1-x_2}
E^+_{\lambda_1+\nu,\lambda_2+\nu,\dots,\lambda_n+\nu}(x)=
E^+_\lambda(x).
\end{equation}

It is seen from the expression \eqref{dett} for an antideterminant
${\det}^+$ that its expression does not change under permutations of
rows or under permutations of columns. This means that for any permutation
$w\in S_n$ we have
\begin{equation}  \label{permu}
E^+_\lambda(wx)=E^+_\lambda(x),\ \ \ \
E^+_{w\lambda}(x)=E^+_\lambda(x).
\end{equation}
Therefore, it is enough to consider only
symmetric exponential functions $E^+_\lambda(x)$
with $\lambda=(\lambda_1,\lambda_2,\dots,\lambda_n)$ such that
\[
\lambda_1\ge \lambda_2\ge \cdots \ge \lambda_n.
\]
Such $\lambda$ are called {\it dominant}. The set of all dominant
$\lambda$ is denoted by $D_+$.
Below, considering symmetric exponential functions
$E^+_\lambda(x)$, we assume that $\lambda\in D_+$.

Antisymmetric multivariate exponential functions of $x=(x_1,x_2,\dots,x_n)$ are
defined as the functions
\begin{alignat}{2}  \label{det-A}
E^-_{\lambda}(x)\equiv &
E^-_{(\lambda_{1},\lambda_2,\dots ,\lambda_{n})}(x):=
\det \left( e^{2\pi {\rm i}\lambda_ix_j}\right)_{i,j=1}^n  \notag\\
\equiv &
 \det \left(
 \begin{array}{cccc}
 e^{2\pi{\rm i}\lambda_1x_1}& e^{2\pi{\rm i}\lambda_1x_2}&\cdots & e^{2\pi{\rm
i}\lambda_1x_{n}}\\
  e^{2\pi{\rm i}\lambda_2x_1}& e^{2\pi{\rm i}\lambda_2x_2}&\cdots & e^{2\pi{\rm
i}\lambda_2x_{n}}\\
  \cdots & \cdots &
\cdots & \cdots \\
 e^{2\pi{\rm i}\lambda_{n}x_1}& e^{2\pi{\rm i}\lambda_{n}x_2}&\cdots & e^{2\pi{\rm
i}\lambda_{n}x_{n}} \end{array} \right)
\notag\\
\equiv & \sum_{w\in S_n} (\det w)
 e^{2\pi{\rm i}\lambda_1x_{w(1)}} e^{2\pi{\rm i}\lambda_2x_{w(2)}}\cdots  e^{2\pi{\rm
i}\lambda_nx_{w(n)}}  
= \sum_{w\in S_n} (\det w) e^{2\pi{\rm i}\langle \lambda, wx \rangle} ,
 \end{alignat}
where $\lambda=(\lambda_1,\lambda_2,\dots,\lambda_n)$ is a set of real
numbers. When $\lambda_1,\lambda_2,\dots,\lambda_n$ are
integers, we denote them as $m=(m_1,m_2,\dots,m_n)$,
\begin{equation}  \label{detII}
E_m^-(x)=\det \left( e^{2\pi{\rm
i}m_ix_j}\right)_{i,j=1}^{n} .
 \end{equation}

It is seen from properties of determinants
that the functions $E^-_\lambda(x)$ satisfy the relation
\begin{equation}  \label{x_1+x_2+}
E^-_\lambda(x_1+a,x_2+a,\dots,x_n+a)=
e^{2\pi{\rm i}(\lambda_1+\lambda_2+\cdots +\lambda_n)a}
E^-_\lambda(x),
\end{equation}
that is, it is enough to consider functions
$E^-_\lambda(x)$ on some hyperplane
$\mathcal{H}_b$. As in the case of symmetric exponential
functions, a transition from one hyperplane
$\mathcal{H}_b$ to another $\mathcal{H}_c$ for the function
$E^-_\lambda$ is fulfilled by means of
multiplication by a usual exponential function
$e^{2\pi{\rm i}\vert \lambda \vert (c-b)}$,
where $\vert \lambda \vert=\lambda_1+\lambda_2+\cdots +\lambda_n$.
For points $x$ of the hyperplane $\mathcal{H}_0$ we have the relation
\begin{equation}  \label{x_1-x_2+}
E^-_{\lambda_1+\nu,\lambda_2+\nu,\dots,\lambda_n+\nu}(x)=
E^-_\lambda(x).
\end{equation}

It follows from properties of determinants that
$E^-_\lambda(x)=0$
if $\lambda$ has at least two coinciding numbers or if $x$
has at least two coinciding coordinates.
For any permutation $w\in S_n$ we receive
\begin{equation}  \label{permuu}
E^-_{w\lambda}(x)=(\det w) E^-_\lambda(x), \ \ \ \
E^-_{\lambda}(wx)=(\det w)E^-_\lambda(x).
\end{equation}
This means that it is enough to consider
antisymmetric exponential functions $E^-_\lambda(x)$
for $\lambda=(\lambda_1,\lambda_2,\dots,\lambda_n)$ such that
\[
\lambda_1> \lambda_2> \cdots > \lambda_n.
\]
Such $\lambda$ are called {\it strictly dominant}. The set of these
$\lambda$ is denoted by $D_+^+$.

\subsection{Affine symmetric group and fundamental domains}
We have seen that the functions $E^+_\lambda(x)$
are symmetric with respect to the permutation group $S_n$, that is,
$E^+_\lambda(wx)=E^+_\lambda(x)$,
$w\in S_n$.
The symmetric exponential functions
$E^+_m$ with integral $m=(m_1,m_2,\dots, m_n)$ admit
additional symmetries related to the periodicity of the exponential functions
$e^{2\pi{\rm i}r y}$, $r\in \mathbb{Z}$, $y\in \mathbb{R}$. These symmetries
are described by the discrete group of shifts in the space $E_n$ by vectors
\[
r_1{\bf e}_1+r_2{\bf e}_2+\cdots +r_n{\bf e}_n,\ \ \ r_i\in \mathbb{Z},
\]
where ${\bf e},{\bf e}_2,\dots, {\bf e}_n$ are the unit vectors in
directions of the corresponding axes. We denote this group by $T_n$.
Permutations of $S_n$ and shifts of $T_n$ generate a group which is denoted
as $S_n^{\rm aff}$ and is called the {\it affine symmetric group}. The
group $S_n^{\rm aff}$ is a semidirect product of its subgroups $S_n$ and $T_n$,
\[
S_n^{\rm aff}=S_n\times T_n,
\]
where $T_n$ is an invariant subgroup, that is, $wtw^{-1}\in T_n$ for $w\in S_n$ and
$t\in T_n$.

An open connected simply connected set $F\subset \mathbb{R}^n$ is called a
{\it fundamental domain} for the group $S_n^{\rm aff}$ (for the group
$S_n$) if it does not contain equivalent points (that is, points $x$
and $x'$ such that $x'=wx$, where $w$ belongs to $S_n^{\rm aff}$ or $S_n$,
respectively) and if its closure contains at least one
point from each $S_n^{\rm aff}$-orbit (from each $S_n$-orbit). Recall that
a $S_n^{\rm aff}$-orbit of a point $x\in \mathbb{R}^n$ is the set of
points $wx$, $w\in S_n^{\rm aff}$.

It is evident that the set $D_+^+$ of all points $x=(x_1,x_2,\dots,x_n)$ such that
\[
x_1>x_2>\cdots >x_n
\]
constitute a fundamental domain for the group $S_n$ (we denote it as
$F(S_n)$). The set of points
$x=(x_1,x_2,\dots,x_n)\in D^+_+$ such that
\[
1>x_1>x_2>\cdots >x_n>0
\]
constitute a fundamental domain for the affine group $S^{\rm aff}_n$
(we denote it as $F(S^{\rm aff}_n)$).

As we have seen, the functions
$E^+_\lambda(x)$
are symmetric with respect to the permutation group $S_n$. This means that
it is enough to consider the functions $E^+_\lambda(x)$
only on the closure of the fundamendal domain $F(S_n)$. Values of
$E^+_\lambda$ on other points are received by using
symmetricity.

Symmetricity of functions $E^+_m$ with integral
$m=(m_1,m_2,\dots,m_n)$ with respect to the affine symmetric group
$S_n^{\rm aff}$,
\begin{equation}\label{aff-sym-aff}
E^+_m(wx+r)=E^+_m(x),\ \ \ \
w\in S_n,\ \ \ r\in T_n,
\end{equation}
means that we may consider
$E^+_m(x)$ only on the closure
of the fundamental domain $F(S^{\rm aff}_n)$, that is, on the set of points
$x$ such that
$1\ge x_1\ge x_2\ge \cdots \ge x_n\ge 0$.
Values of
$E^+_m(x)$ on other points are obtained by
using the relation \eqref{aff-sym-aff}.

The exponential functions
$E^-_m(x)$ with integral $m=(m_1,m_2,\dots, m_n)$ also admit
additional symmetries related to the periodicity of the usual exponential functions
$e^{2\pi{\rm i}r y}$, $r\in \mathbb{Z}$, $y\in \mathbb{R}$. These symmetries
are described by the affine symmetric group $S^{\rm aff}_n$.
We have
\begin{equation}\label{aff-sym-af}
E^-_m(wx+r)=(\det w)E^-_m(x),\ \ \ \
w\in S_n,\ \ \ r\in T_n,
\end{equation}
that is, it is enough to consider the functions
$E^-_m(x)$ only on the closure of the fundamental
domain $F(S^{\rm aff}_n)$. Values of $E^-_m(x)$
on other points are obtained by using the relation
\eqref{aff-sym-af}.

The functions
$E^-_\lambda(x)$, $\lambda=(\lambda_1,\lambda_2,\dots, \lambda_n)$,
$\lambda_i\in \mathbb{R}$, are antisymmetric with respect to
the symmetric group $S_n$,
\[
E^-_\lambda(wx)=(\det w) E^-_m(x),\ \ \ \
w\in S_n.
\]
For this reason, we may consider $E^-_\lambda$
only on the fundamental domain $F(S_n)$.

\subsection{Properties}
Symmetricity and antisymmetricity of symmetric and antisymmetric
multivariate exponential
functions are  main properties of these functions. However, they possess
many other interesting properties.
 \medskip

{\bf Behavior on boundary.}
The symmetric and antisymmetric functions $E^+_\lambda(x)$
and $E^-_\lambda(x)$
are  finite sums of exponential functions. Therefore, they are continuous
functions of $x_1,x_2,\dots,x_n$ and have
continuous derivatives of all orders in $\mathbb{R}^n$.

The closure $\overline{F(S_n)}$ of the fundamental domain $F(S_n)$
without points of $F(S_n)$ is called
a {\it boundary} of the fundamental domain $F(S_n)$ and is
denoted by $\partial F(S_n)$. A point $x=(x_1,x_2,\dots, x_n)\in
\overline{F(S_n)}$ belongs to $\partial F(S_n)$ if and only if
at least two coordinates $x_i,x_j$ in $x$ coincide. It is clear that
the boundary $\partial F(S_n)$ is composed of points of
$\overline{F(S_n)}$
belonging to the hyperplanes given by the equations
\[
x_i=x_j,\ \ \ \ i,j=1,2,\dots,n,\ \ \ i\ne j.
\]

Similarly, the boundary $\partial F(S^{\rm aff}_n)$ of the
fundamental domain $F(S^{\rm aff}_n)$ consists of points
of $\overline{F(S^{\rm aff}_n)}$ which do not belong to
$F(S^{\rm aff}_n)$. A point
$x=(x_1,x_2,\dots, x_n)\in
\overline{F(S^{\rm aff}_n)}$ belongs to $\partial F(S^{\rm aff}_n)$
if and only if at least two coordinates $x_i,x_j$ in $x$ coincide
or if one of the conditions $x_1=1$, $x_n=0$ is fulfilled.

It follows from properties of determinants that the function
$E^-_\lambda(x)$ vanishes on the boundary
$\partial F(S_n)$,
\[
E^-_\lambda(x)=0,\ \ \ {\rm for}\ \ \ x\in \partial F(S_n).
\]
This relation is true for $E^-_m(x)$, $m_i\in \mathbb{Z}$.
In this case we also have $E^-_m(x)=0$ for points $x\in \partial F(S^{\rm aff}_n)$
such that $x_1-x_n=1$.

For symmetric multivariate functions $E^+_\lambda(x)$ we have
\[
 \frac{\partial E^+_\lambda(x)}{\partial {\bf n}}=0\ \ \
{\rm for}\ \ \ x\in \partial F(S_n)
\]
where ${\bf n}$ is the normal to the boundary
$\partial F(S_n)$.

 \medskip

{\bf Complex conjugation.}
Let $\lambda=(\lambda_1,\lambda_2,\dots ,\lambda_{n})$ be a strictly dominant
element, that is, $\lambda_1> \lambda_2> \cdots > \lambda_{n}$.
We have
\begin{gather}\label{compl-con}
E^-_\lambda(x)  =
 \sum _{w\in S_{n}} (\det w) e^{2\pi{\rm i}((w\lambda)_1x_1+
 \cdots +(w\lambda)_{n}x_{n})} ,
\end{gather}
where $(w\lambda)_{1},(w\lambda)_{2},\dots ,(w\lambda)_{n}$
are the coordinates of the point $w\lambda$.

The element $-(\lambda_{n},\lambda_{n-1},\dots ,\lambda_1)$ is strictly
dominant if the element $(\lambda_1,\lambda_2,\dots$, $\lambda_{n})$ is strictly
dominant. In the group $S_n$ there exists an element $w_0$
such that
\[
w_0(\lambda_{1},\lambda_2,\dots ,\lambda_{n})=
(\lambda_{n},\lambda_{n-1},\dots ,\lambda_1)
\]
It is easy to calculate that the set
$(\lambda_{n},\lambda_{n-1},\dots ,\lambda_1)$ is obtained from
$(\lambda_{1},\lambda_2,\dots ,\lambda_{n})$ by $n(n-1)/2$ permutations
of two neighboring numbers. Clearly, $\det w_0=1$ if this number
is even and $\det w_0=-1$ otherwise. Thus,
\begin{gather*}
\det w_0=1\qquad {\rm for}\qquad n=4k\qquad {\rm and}\qquad
n=4k+1,
\\
\det w_0=-1\qquad {\rm for}\qquad n=4k-2\qquad {\rm and}\qquad
n=4k-1,
\end{gather*}
where $k$ is an integer.
It follows from here that in the expressions for the
exponential functions $E^-_{(\lambda_{1},\lambda_2,
\dots ,\lambda_{n})}(x)$ and
$E^-_{-(\lambda_{n},\lambda_{n-1},\dots ,
\lambda_{1})}(x)$ there are summands
\begin{gather}\label{w_0}
e^{2\pi{\rm i}\langle w_0\lambda ,x\rangle}= e^{2\pi{\rm
i}(\lambda_{n}x_1+
 \cdots +\lambda_{1}x_{n})}\qquad
{\rm and}\qquad e^{-2\pi{\rm i}(\lambda_{n}x_1+
 \cdots +\lambda_{1}x_{n})} ,
\end{gather}
respectively, which are complex conjugate to each other. Moreover,
the first expression is contained with the sign $(\det w_0)$ in
$E^-_{(\lambda_{1},\lambda_2,\dots ,
\lambda_{n})}(x)$, that is, the
expressions \eqref{w_0} are contained in
$E^-_{(\lambda_{1},\lambda_2,\dots ,\lambda_{n})}(x)$ and
$E^-_{-(\lambda_{n},\lambda_{n-1},
\dots ,\lambda_{1})}(x)$ with the same sign for
$n=4k, 4k+1$ and with opposite signs for $n=4k-2, 4k-1$, $k\in
\mathbb{Z}$.

 Similarly, in the expressions \eqref{compl-con} for the function
$E^-_{(\lambda_{1},\lambda_2,\dots ,\lambda_{n})}(x)$
and for the function $E^-_{-(\lambda_{n},\lambda_{n-1},
\dots ,\lambda_{1})}(x)$ all other summands are
(up to a sign, which depends on a value of $n$) pairwise complex
conjugate. Therefore,
\begin{gather}\label{compl}
E^-_{(\lambda_{1},\lambda_2,\dots
,\lambda_{n})}(x)=\overline{E^-_{-(\lambda_{n},
\lambda_{n-1},\dots ,\lambda_{1})}(x)}
\end{gather}
for $n=4k,4k+1$ and
\begin{gather}\label{compl-n}
E^-_{(\lambda_{1},\lambda_2,\dots
,\lambda_{n})}(x)=-\overline{E^-_{-(\lambda_{n},
\lambda_{n-1},\dots ,\lambda_{1})}(x)}
\end{gather}
for $n=4k-2,4k-1$.

According to \eqref{compl} and \eqref{compl-n}, if
\begin{gather}\label{m-A}
(\lambda_1,\lambda_2,\dots ,\lambda_{n})=-(\lambda_{n},\lambda_{n-1},
\dots ,\lambda_{1}),
\end{gather}
then {\it the  function $E^-_\lambda$ is real for
$n=4k,4k+1$ and pure imaginary for $n=4k-2,4k-1$.}
Moreover, the right hand side of \eqref{compl-con} for this case consists of
pairs of terms which give sine or cosine functions.
It is representable as a sum
of cosines of angles if $n=4k,4k+1$ and
as a sum of sines of
angles multiplied by ${\rm i}=\sqrt{-1}$ if $n=4k-2,4k-1$.

It is proved similarly that for the symmetric exponential functions
$E^+_\lambda(x)$ and $E^+_{w_0\lambda}(x)$
we have the following relation
\[
E^+_{(\lambda_{1},\lambda_2,\dots
,\lambda_{n})}(x)=\overline{E^+_{-(\lambda_{n},
\lambda_{n-1},\dots ,\lambda_{1})}(x)}.
\]
If
$(\lambda_1,\lambda_2,\dots ,\lambda_{n})=-(\lambda_{n},\lambda_{n-1},
\dots ,\lambda_{1})$, then the function
$E^+_{\lambda}(x)$ is real. In this case
$E^+_{\lambda}(x)$ can be represented as a sum of
cosines of the corresponding angles.
 \medskip

{\bf Scaling symmetry.}
For $c\in \mathbb{R}$, let $c\lambda=(c\lambda_1,c\lambda_2,\dots,c\lambda_n)$.
Then
\[
E^-_{c\lambda}(x)=\sum_{w\in W} (\det w) e^{2\pi{\rm i}\langle
cw\lambda ,x\rangle} = \sum_{w\in W} (\det w) e^{2\pi{\rm
i}\langle w\lambda ,cx\rangle} = E^-_{\lambda}(cx).
\]
The equality $E^-_{c\lambda}(x)=E^-_{\lambda}(cx)$
expresses the {\it scaling symmetry of exponential functions
$E^-_{c\lambda}(x)$.} The scaling symmetry is true also for symmetric
exponential functions,
$E^+_{c\lambda}(x)=E^+_{\lambda}(cx)$.
 \medskip

{\bf Duality.}
Due to invariance of the scalar product $\langle \cdot, \cdot
\rangle$ with respect to the symmetric group $S_n$, $\langle w\mu ,wy
\rangle= \langle \mu ,y \rangle$, for $x=(x_1,x_2,\dots,x_n)$, $x_i\ne x_j$,
$i\ne j$, we have
 \[
E^-_\lambda(x)=\sum_{w\in W}(\det w) e^{2\pi {\rm i}\langle
\lambda,w^{-1}x \rangle} =\sum_{w\in W}(\det w) e^{2\pi {\rm
i}\langle x, w\lambda \rangle}=  E^-_x(\lambda).
 \]
This relation expresses
the {\it duality} of antisymmetric orbit functions.
The duality is true also for symmetric exponential functions,
$E^+_{\lambda}(x)=E^+_{x}(\lambda)$.
 \medskip

{\bf Orthogonality on the fundamental domain $F(S_n^{\rm aff})$.}
Antisymmetric exponential functions $E^-_m$ with
$m=(m_1,m_2,\dots,m_n)\in D^+_+$, $m_j\in \mathbb{Z}$,
are orthogonal on $F(S_n^{\rm aff})$ with respect to the Euclidean
measure,
 \begin{gather}\label{ortogo}
|F(S_n^{\rm aff})|^{-1}\int_{F(S_n^{\rm aff})}
E^-_m(x)\overline{ E^-_{m'}(x)}dx=
|S_n|\delta_{mm'} ,
 \end{gather}
where the overbar means complex conjugation, $|S_n|$ means a number of
elements in the set $S_n$, and $|F(S_n^{\rm aff})|$ is an area of the
fundamental domain $F(S_n^{\rm aff})$.
This relation follows from the equality
\[
\int_{\sf T} E^-_m(x)
\overline{E^-_{m'}(x)}dx=
|S_n|\delta_{mm'}
 \]
(where ${\sf T}$ is the torus in $E_n$ consisting of points
$x=(x_1,x_2,\dots,x_n)$, $0\le x_i<1$), which is a consequence of
orthogonality of the exponential functions
$e^{2\pi{\rm i}\langle \mu,x \rangle}$ (entering into the
definition of $E^-_m(x)$) for different sets $\mu$.

If to assume that an area of
${\sf T}$ is equal to 1, $|{\sf T}|=1$, then $|F(S_n^{\rm aff})|=|S_n|^{-1}$ and
formula \eqref{ortogo} takes the form
\begin{gather}\label{ortog-tor}
\int_{F(S_n^{\rm aff})} E^-_m(x)
\overline{E^-_{m'}(x)}dx=
\delta_{mm'} .
 \end{gather}

In the expression \eqref{det-SS} for symmetric exponential functions
there can be coinciding summands. For this reason, for
symmetric exponential functions
 $E^+_m$ with
$m=(m_1,m_2,\dots,m_n)\in D_+$, $m_j\in \mathbb{Z}$,
the relation \eqref{ortog-tor} is replaced by
\begin{gather}\label{ortog-torus}
\int_{F(S_n^{\rm aff})} E^+_m(x)
\overline{E^+_{m'}(x)}dx=
|S_m| \delta_{mm'} .
 \end{gather}
where $|S_m|$ is a number of elements in the subgroup $S_m$ of
$S_n$ consisting of elements $w\in S_n$ such that $wm=m$.
 \medskip

{\bf Orthogonality of symmetric and antisymmetric exponential
functions.}
Let $w_i$ ($i=1,2,\dots,n-1$) be the permutation of coordinates $x_i$ and
$x_{i+1}$. We create the domain $F^{\rm
ext}(S_n^{\rm aff})=F(S_n^{\rm aff})\cup w_i F(S_n^{\rm aff})$,
where $F(S_n^{\rm aff})$ is the
fundamental domain for the affine group $S_n^{\rm aff}$. Since
for $m=(m_1,m_2,\dots,m_n)\in D^+_+$, $m_j\in \mathbb{Z}$, we have
$E^+_m(w_i x)=E^+_m(x)$ and $E^-_m(w_i x)=-E^-_m(x)$,
then
 \begin{gather}\label{ortog-ant}
 \int_{F^{\rm ext}(S_n^{\rm aff})} E^+_m(x)
 \overline{E^-_{m'}(x)}dx=0.
 \end{gather}
Indeed, due to symmetry and antisymmetry of symmetric and
antisymmetric exponential functions, respectively, we have
\[
\int_{F^{\rm ext}(S_n^{\rm aff})} E^+_m(x)
\overline{E^-_{m'}(x)}dx
= \int_{F(S_n^{\rm aff})} E^+_m(x)
\overline{E^-_{m'}(x)}dx+
\int_{w_i F(S_n^{\rm aff})} E^+_m(x)
\overline{E^-_{m'}(x)}dx
\] \[   \qquad\qquad\qquad
= \int_{F(S_n^{\rm aff})}
 E^+_m(x)
\overline{E^-_{m'}(x)}dx+ \int_{F(S_n^{\rm aff})}
E^+_m(x)\overline{(-E^-_{m'}(x)})dx=0.
\]
The relation \eqref{ortog-ant} is a generalization of
the orthogonality of the functions sine and cosine on the interval
$(0,2\pi)$.

\subsection{Special cases}
The special case of symmetric and antisymmetric exponential
functions at $(\lambda_1,\lambda_2,\dots,\lambda_n)=
\frac12 (n-1,n-3,\dots,-n+3,-n+1)\equiv \rho$
is of great interest since it is met in the representatin theory.
The antisymmetric exponential function $E^-_\rho(x)$ is given by the formula
\begin{equation} \label{a---A}
E^-_\rho(x)=(2{\rm i})^{n(n-1)/2} \prod_{1\le i<j\le n} \sin\,
\pi (x_i-x_j).
 \end{equation}
It follows if to represent  $\sin\, \pi (x_i-x_j)$ in terms of
exponential functions, then to fulfil multiplication of these functions
and to compare with the expression \eqref{det-A} for $E^-_\rho(x)$.

Let us set $(m_1,m_2,\dots,m_n)=(n-1,n-2,\dots,1,0)\equiv \rho'$.
The antisymmetric exponential function $E^-_{\rho'}(x)$
can be written down  in the form of the
Vandermonde determinant,
 \begin{equation} \label{a-pho-A}
E^-_{\rho'}(x)=\det \left( e^{2\pi{\rm i}(n-i)x_j}
\right)_{i,j=1}^{n}
 =\prod_{k<l} (e^{2\pi{\rm i}x_k}-e^{-2\pi{\rm i}x_l}).
 \end{equation}
The last equality follows from the expression for the
Vandermonde determinant. Since $\rho'=\rho+\frac{n-1}{2}$, the expressions
\eqref{a---A} and \eqref{a-pho-A} are connected by the relation
\[
 E^-_{\rho'}(x)=e^{\pi {\rm i}|x| (n-1)}E^-_{\rho}(x),
\]
where $|x|=\sum_{i=1}^n x_i$.

It is easy to see that {\it the function $E^-_\rho(x)$
does not vanish on intrinsic points of the fundamental domain $F(S^{\rm aff}_n)$.}

The symmetric counterpart $E^+_\rho(x)$ of the formula \eqref{a---A} for
the antisymmetric exponential function $E^-_\rho(x)$ has the form
\begin{equation}\label{rho-symm}
E^+_\rho(x)=2^{n(n-1)/2} \prod_{1\le i<j\le n} \cos \, \pi
(x_i-x_j).
\end{equation}

\section{Solutions of the Laplace equation}
The Laplace operator on the $n$-dimensional Euclidean space
$E_n$ in the orthogonal coordinates
$x=(x_1,x_2,\dots,x_n)$ has the
form
\[
\Delta=\frac{\partial^2}{\partial
x^2_1}+\frac{\partial^2}{\partial x^2_2}+\cdots
+\frac{\partial^2}{\partial x^2_n} .
\]

We take any summand in the expression for symmetric or antisymmetric
multivariate exponential function
and act upon it by the operator
$\Delta$. We get
\begin{gather*}
\Delta  e^{2\pi{\rm i}((w( \lambda))_1x_1+
 \cdots + (w(\lambda))_nx_{n})}
=-4\pi^2 \langle \lambda,\lambda \rangle\, e^{2\pi{\rm
i}((w( \lambda))_1x_1+
 \cdots + (w(\lambda))_nx_{n})},
\end{gather*}
where $\lambda=(\lambda_1,\lambda_2,\dots ,\lambda_n)$  determines
$E^-_\lambda(x)$.
Since this action of $\Delta$ does not depend on a summand
from the expression for symmetric or antisymmetric
exponential function, we have
\begin{gather}\label{Lapl}
\Delta E^-_\lambda(x)= -4\pi^2\langle \lambda,\lambda \rangle
 E^-_\lambda(x),\ \ \ \
\Delta E^+_\lambda(x)= -4\pi^2\langle \lambda,\lambda \rangle
 E^+_\lambda(x).
 \end{gather}

The formula \eqref{Lapl} can be generalized in the following way.
Let $\sigma_k(y_1,y_2,\dots,y_n)$ be the $k$-th elementary
symmetric polynomial of degree $k$, that is,
\[
 \sigma_k(y_1,y_2,\dots,y_n)=\sum_{1\le k_1<k_2<\cdots <k_n\le n}
y_{k_1} y_{k_2}\cdots y_{k_n}.
\]
Then for $k=1,2,\dots,n$ we have
\begin{gather}\label{Lap-gene}
\sigma_k\left( \tfrac{\partial^2}{\partial
x^2_1},\tfrac{\partial^2}{\partial
x^2_2},\dots,\tfrac{\partial^2}{\partial x^2_n} \right)
E^{\pm}_\lambda(x)= (-4\pi^2)^k \sigma_k
(\lambda_1^2,\lambda_2^2,\dots,\lambda_n^2)
E^{\pm}_\lambda(x).
 \end{gather}
Note that $n$ differential equations \eqref{Lap-gene} are algebraically
independent.

Thus, antisymmetric exponential functions
$E^-_\lambda(x)$
are eigenfunctions of the
operators $\sigma_k\left( \frac{\partial^2}{\partial
x^2_1},\frac{\partial^2}{\partial
x^2_2},\dots,\frac{\partial^2}{\partial x^2_n} \right)$,
$k=1,2,\dots,n$, on the fundamental domain $F(S_n)$ of the symmetric group
$S_n$ satisfying the boundary condition
\begin{gather}\label{Neum1}
E^-_m(x) =0 \qquad {\rm for}\qquad x\in \partial F(S_n).
\end{gather}

Similarly, symmetric exponential functions
$E^+_\lambda(x)$ are eigenfunctions of the
operators $\sigma_k\left( \frac{\partial^2}{\partial
x^2_1},\frac{\partial^2}{\partial
x^2_2},\dots,\frac{\partial^2}{\partial x^2_n} \right)$,
$k=1,2,\dots,n$, on the fundamental domain $F(S_n)$ satisfying the
boundary condition
$$
\frac{\partial\, E^+_m(x)}{\partial\, {\bf n}}
 =0 \qquad {\rm for}\qquad x\in \partial F(S_n),
$$
where ${\bf n}$ is the normal to the boundary $\partial F(S_n)$.
That is, they are solutions of the Neumann boundary value problem for the
domain $F(S_n)$.

\section{Symmetric and antisymmetric Fourier transforms}\label{section10}

Symmetric and antisymmetric exponential functions determine
symmetric and antisymmetric multivariate Fourier transforms which
generalize the usual Fouirier transform.

As in the case of exponential functions of one variable, (anti)symmetric
exponential functions determine three types of Fourier transforms:
 \medskip

(a) Fourier transforms related to the exponential
functions $E^{\pm}_m(x)$ with $m=(m_1,m_2,\dots,m_n)$,
$m_j\in \mathbb{Z}$ (Fourier series);
 \medskip

(b) Fourier transforms related to $E^{\pm}_\lambda(x)$
with $\lambda\in D_+$;
 \medskip

(c) symmetrized and antisymmetrized multivariate finite Fourier transforms.

\subsection{Expansions in (anti)symmetric
exponential functions on $F(S^{\rm aff}_n)$}
Let $f(x)$ be symmetric (with respect to the affine
symmetric group $S_n^{\rm aff}$) continuous function on the
$n$-dimensional Euclidean space $E_n$
which has continuous derivatives. We may consider this
function on the torus ${\sf T}$ which is a closure
of the union of the sets $w F(S^{\rm aff}_n)$, $w\in S_n$.
The function $f(x)$, as a function on ${\sf T}$, can be expanded
in exponential functions
$e^{2\pi{\rm i}m_1 x_1}
e^{2\pi{\rm i}m_2 x_2}\cdots e^{2\pi{\rm i}m_n x_n},\ \ \
m_i\in \mathbb{Z}$.
We have
 \begin{equation}\label{expan-anti}
f(x)=\sum_{m_i\in \mathbb{Z}} c_m
e^{2\pi{\rm i}m_1 x_1}
e^{2\pi{\rm i}m_2 x_2}\cdots e^{2\pi{\rm i}m_n x_n},
 \end{equation}
where $m=(m_1,m_2,\dots,m_n)$.
Due to symmetry $f(wx)=f(x)$, $w\in S_n$, for any $w\in S_n$ we have
\[
f(wx)=
\sum_{m_i\in \mathbb{Z}} c_m
e^{2\pi{\rm i}m_1 x_{w(1)}}
\cdots e^{2\pi{\rm i}m_n x_{w(n)}}
=\sum_{m_i\in \mathbb{Z}} c_{m}
e^{2\pi{\rm i}m_{w^{-1}(1)} x_{1}}
\cdots
e^{2\pi{\rm i}m_{w^{-1}(n)} x_n}
\] \[
=\sum_{m_i\in \mathbb{Z}} c_{wm}
e^{2\pi{\rm i}m_1 x_{1}}
\cdots e^{2\pi{\rm i}m_n x_{n}}
=f(x)=
\sum_{m_i\in \mathbb{Z}} c_m
e^{2\pi{\rm i}m_1 x_1}
\cdots e^{2\pi{\rm i}m_n x_n}.
\]
Therefore, the coefficients $c_m$ satisfy the conditions
$c_{wm}=c_m$, $w\in S_n$.
Collecting in \eqref{expan-anti} exponential functions at the same
$c_{wm}$, $w\in S_n$, we obtain the expansion
 \begin{equation}\label{expan-anti-2}
f(x)=\sum_{m\in P_+} c_m
E^+_m(x),
 \end{equation}
where $P_+=D_+\cap \mathbb{Z}$.
Thus, {\it any symmetric
(with respect to $S_n$) continuous function $f$ on ${\sf T}$
which has continuous derivatives (that is, any continuous function on $D_+$
with continuous derivatives) can be expanded in symmetric
exponential functions} $E^+_m(x)$, $m\in P_+$.

By the orthogonality relation \eqref{ortog-torus}, the coefficients $c_m$ in the
expansion \eqref{expan-anti-2} are determined by the formula
 \begin{gather}\label{decom-ant-2}
c_m =|S_m|^{-1} \int_{F(S^{\rm aff}_n)} f(x)
\overline{E^+_m(x)}dx ,
 \end{gather}
where, as before, $|S_m|$ is a number of elements in the subgroup $S_m$
of $S_n$ consisting of $w\in S_n$ such that $wm=m$.
Moreover, the Plancherel formula
 \begin{gather}\label{decom-ant-3}
\sum_{m\in P^+}  |c_m|^2=|S_m|^{-1}  \int_{F(S^{\rm aff}_n)} |f(x)|^2dx
 \end{gather}
holds, which means that the Hilbert spaces with the appropriate
scalar products are isometric.

Formula \eqref{decom-ant-2} is the symmetrized Fourier transform
of the function $f(x)$. Formula \eqref{expan-anti-2} gives an inverse
transform. Formulas \eqref{expan-anti-2} and \eqref{decom-ant-2} give the
{\it symmetric multivariate Fourier transforms} corresponding to the symmetric
exponential functions $E^+_m(x)$, $m\in P^+$.

Analogous transforms hold for antisymmetric
exponential functions $E^-_m(x)$, $m\in P^+_+\equiv D^+_+\cap \mathbb{Z}$.
Let $f(x)$ be antisymmetric (with respect to the symmetric
group $S_n$) continuous function on the
$n$-dimensional torus ${\sf T}$,
which has continuous derivatives. We may consider this function as a
function on $F(S^{\rm aff}_n)$. Then we have the expansion
 \begin{equation}\label{expan-anti-sym-2}
f(x)=\sum_{m\in P_+^+} c_m
E^-_m(x),
\ \ \
{\rm where}\ \ \
c_m = \int_{F(S^{\rm aff}_n)} f(x)
\overline{E^-_m(x)}dx.
 \end{equation}
Moreover, the Plancherel formula holds:
 \begin{gather}\label{decom-ant-sym-3}
\sum_{m\in P^+_+}  |c_m|^2= \int_{F(S^{\rm aff}_n)} |f(x)|^2dx.
 \end{gather}

Let ${\mathcal L}_0^2(F(S^{\rm aff}_n))$ denote the Hilbert space of functions on
the fundamental domain $F(S^{\rm aff}_n)$, which behave on the boundary $\partial
F(S^{\rm aff}_n)$ of the fundamental domain $F(S^{\rm aff}_n)$ in the
same way as the functions $E^-_m(x)$ do. Let
\[
\langle f_1,f_2\rangle = \int_{F(S^{\rm aff}_n)} f_1(x)\overline{f_2(x)} dx
\]
be a scalar product in this space. The formulas
\eqref{expan-anti-sym-2}-\eqref{decom-ant-sym-3} show that {\it the set of
exponential
functions $E^-_m(x)$, $m\in P^+_+$, form an
orthogonal basis of ${\mathcal L}_0^2(F(S^{\rm aff}_n))$.}

Let $F^{\rm ext}(S^{\rm aff}_n)=F(S^{\rm aff}_n)\bigcup F(w_i S^{\rm aff}_n)$
denote the set from section 2. Then we can extend the
symmetric and antisymmetric Fourier transforms to the functions
from the Hilbert space
${\mathcal L}^2(F^{\rm ext}(S^{\rm aff}_n))$  with the scalar product
\[
\langle f_1,f_2\rangle = \int_{F^{\rm ext}(S^{\rm aff}_n)}
f_1(x)\overline{f_2(x)} dx.
\]
This transform is of the form
 \begin{equation}\label{expan-a-s-2}
f(x)=\sum_{m\in P_+} c_m
E^+_m(x)+ \sum_{m\in P^+_+} c'_m
E^-_m(x),
 \end{equation}
where
 \begin{gather}\label{decom-an-s-2}
c_m =|S_m|^{-1} \int_{F(S^{\rm aff}_n)} f(x)
\overline{E^+_m(x)}dx , \ \ \
c'_m =\int_{F(S^{\rm aff}_n)} f(x)
\overline{E^-_m(x)}dx.
 \end{gather}
The corresponding Plancherel formula holds.
The functions $E^+_m(x)$, $m\in P_+$, and $E^-_m(x)$, $m\in P^+_+$,
form a complete orthogonal basis of the Hilbert space
${\mathcal L}^2(F^{\rm ext}(S^{\rm aff}_n))$.

\subsection{Multivariate Fourier transforms on the fundamental domain
$F(S_n)$}\label{section10.3}
The expansions \eqref{expan-anti-2} and \eqref{expan-anti-sym-2}
of functions on the fundamental
domain $F(S^{\rm aff}_n)$ are respectively expansions in the symmetric
and antisymmetric exponential functions $E^+_m(x)$ and
$E^-_m(x)$ with integral
$m=(m_1,m_2,\dots,m_n)$.  The exponential functions
$E^+_\lambda(x)$ and
$E^-_\lambda(x)$ with $\lambda$ lying in the fundamental
domain $F(S_n)$ (and not
obligatory integral) are not invariant (anti-invariant) with respect to the
corresponding affine symmetric group $S_n^{\rm aff}$. They are
invariant (anti-invariant) only with respect to the permutation group
$S_n$. A fundamental
domain of $S_n$ coincides with the set $D^+_+$ consisting of the points $x$ such that
$m_1>m_2>\cdots >m_n$. For this reason,
the functions $E^-_\lambda(x)$,
$\lambda\in D^+_+$, and $E^+_\lambda(x)$,
$\lambda\in D_+$,  determine Fourier transforms on $D_+$.

We began with the usual Fourier transforms on ${\mathbb R}^n$:
 \begin{gather}\label{F-1}
\tilde f (\lambda)=\int_{{\mathbb R}^n} f(x) e^{2\pi {\rm
i}\langle \lambda,x \rangle} dx,
  \\
  \label{F-2}
 f (x)=\int_{{\mathbb R}^n} \tilde f(\lambda) e^{-2\pi {\rm i}\langle
\lambda,x \rangle} d\lambda.
  \end{gather}
Let the function $f(x)$ be anti-invariant with respect to the symmetric
group $S_n$, that is, $f(wx)=(\det w)f(x)$, $w\in S_n$. It is easy to
check that $\tilde f (\lambda)$ is also
anti-invariant with respect to the group $S_n$. Replace in
\eqref{F-1} $\lambda$ by $w\lambda$, $w\in S_n$, multiply both sides
by $\det w$, and sum these both side over $w\in S_n$. Then instead
of \eqref{F-1} we obtain
 \begin{gather}\label{F-3}
\tilde f (\lambda)= \int_{D_+} f(x) E^-_\lambda(x) dx,\qquad
\lambda\in D^+_+,
  \end{gather}
where we have taken into account that $f(x)$ is anti-invariant
with respect to $S_n$.

Similarly, starting from \eqref{F-2}, we obtain the inverse
formula:
 \begin{gather}\label{F-4}
 f (x)= \int_{D_+} \tilde f(\lambda)
 \overline{E^-_\lambda(x)} d\lambda .
  \end{gather}
For the transforms \eqref{F-3} and \eqref{F-4} the Plancherel
formula
 \[
 \int_{D_+} |f(x)|^2 dx=
\int_{D_+} |\tilde f(\lambda) |^2  d\lambda
 \]
holds. Formulas \eqref{F-3} and \eqref{F-4} determine the {\it
antisymmetric multivariate Fourier transform on the domain}
$F(S_n)$.

Similarly, starting from formulas
\eqref{F-1} and \eqref{F-2} we receive the symmetric multivariate
Fourier transform on the domain $F(S_n)$:
 \begin{gather}\label{F-5}
\tilde f (\lambda)= \int_{D_+} f(x) E^+_\lambda(x) dx,\qquad
 f (x)= \int_{D_+} \tilde f(\lambda)
 \overline{E^+_\lambda(x)} d\lambda .
  \end{gather}
The corresponding Plancherel formula holds.

\section{Multivariate antisymmetric and symmetric\\
finite Fourier transforms}
Along with the integral Fourier transform in one variable there exists a discrete
Fourier transform in one variable running over a finite set.
Similarly, it is possible to introduce  finite multivariate
antisymmetric and symmetric Fourier transforms, based on
antisymmetric and symmetric exponential functions.
We first consider the finite Fourier transform in one
variable, which will be used below.
Then we expose a general antisymmetric and symmetric Fourier
transforms. Under exposition we use the methods developed in \cite{MP06}.

\subsection{Finite Fourier transform}
Let us fix a positive integer $N$ and consider the numbers
 \begin{equation}\label{f-F-1}
e_{mn}:=N^{-1/2}\exp (2\pi {\rm i}mn/N),\ \ \ \ m,n=1,2,\cdots,N.
 \end{equation}
The matrix $(e_{mn})_{m,n=1}^N$ is unitary, that is,
 \begin{equation}\label{f-F-2}
\sum_k e_{mk}\overline{e_{nk}} =\delta_{mn},\ \ \ \ \sum_k
e_{km}\overline{e_{kn}} =\delta_{mn}.
 \end{equation}
Indeed, according to the formula for a sum of a geometric
progression we have
\[
t^a+t^{a+1}+\cdots +t^{a+r}=(1-t)^{-1}t^a(1-t^{r+1}),\ \ \ \ t\ne
1,
\]  \[
t^a+t^{a+1}+\cdots +t^{a+r}=r+1,\ \ \ \ t=1.
\]
Setting $t=\exp (2\pi{\rm i}(m-n)/N)$, $a=1$ and $r=N-1$, we prove
\eqref{f-F-2}.

Let $f(n)$ be a function of $n\in \{ 1,2\cdots ,N\}$. We may
consider the transform
 \begin{equation}\label{f-F-3}
\sum_{n=1}^N f(n)e_{mn}\equiv N^{-1/2} \sum_{n=1}^N f(n) \exp
(2\pi{\rm i}mn/N) =\tilde f (m).
 \end{equation}
Then due to unitarity of the matrix $(e_{mn})_{m,n=1}^N$, we
express $f(n)$ as a linear combination of conjugates of
the functions \eqref{f-F-1}:
 \begin{equation}\label{f-F-4}
f(n)= N^{-1/2} \sum_{m=1}^N {\tilde f}(m) \exp (-2\pi{\rm i}mn/N)
.
 \end{equation}
The function ${\tilde f}(m)$ is a {\it finite Fourier transform}
of $f(n)$. This transform is a linear map. The
formula \eqref{f-F-4} gives an inverse transform. The Plancherel
formula
\[
\sum_{m=1}^N |\tilde f(m)|^2=\sum_{n=1}^N | f(n)|^2
\]
holds for transforms \eqref{f-F-3} and \eqref{f-F-4}. This means
that the finite Fourier transform conserves the norm introduced in
the space of functions on $\{ 1,2,\dots,N\}$.

\subsection{Antisymmetric multivariate discrete
Fourier transforms}
We use the discrete exponential function \eqref{f-F-1},
\begin{equation}\label{mult-F-1-1-1}
 \textstyle{
e_m(s):=N^{-1/2}\exp (2\pi{\rm i}ms),\ \ \ s\in F_N\equiv\{ \frac1N,
\frac2N,\dots,\frac{N-1}N ,1\},\ \ \ m\in {\mathbb Z}^{\ge 0},
}
\end{equation}
and make a multivariate discrete exponential function by taking
a product of $n$ copies of these functions,
\begin{alignat}{2}\label{mult-F-1}
e_{\bf m}({\bf s}):=&\; e_{m_1}(s_1) e_{m_2}(s_2)\cdots
e_{m_n}(s_n)\notag\\
=&\; N^{-n/2}\exp (2\pi{\rm i}m_1s_1)\exp (2\pi{\rm i}m_2s_2)
\cdots \exp (2\pi{\rm i}m_ns_n)
 \end{alignat}
where ${\bf s}=(s_1,s_2,\dots,s_n)\in F^n_N$ and
${\bf m}=(m_1,m_2,\dots,m_n)\in ({\mathbb Z}^{\ge 0})^n$.
Now we take these multivariate functions for
integers $m_i$ such that
$m_1>m_2>\cdots>m_n\ge 0$ and make an antisymmetrization. As a result, we
obtain a finite version of the antisymmetric exponential functions \eqref{det-A}:
\begin{equation}\label{mult-F-2}
\tilde E^-_{\bf m}({\bf s}):=
  |S_n|^{-1/2} \det (e_{m_i}(s_j))_{i,j=1}^n =|S_n|^{-1/2}N^{-n/2}E^-_{\bf m}({\bf s}),
 \end{equation}
where, as before, $|S_n|$ is the order of the symmetric group $S_n$.

The $n$-tuples ${\bf s}$ in \eqref{mult-F-2} runs over $F_N^n
\equiv F_N\times \cdots \times F_N$ ($n$ times). We
denote by $\hat F_N^n$ the subset of $F_N^n$ consisting of ${\bf s}\in
F_N^n$ such that
$$
s_1>s_2>\cdots >s_n.
$$
The set $\hat F_N^n$ is a finite subset of the fundamental domain
$F(S^{\rm aff}_n)$ of the group $S^{\rm aff}_n$.

Note that acting by permutations $w\in S_n$ upon $\hat F_N^n$ we obtain
the whole set $F_N^n$ without those points which are invariant under
some nontrivial permutation $w\in S_n$. Clearly, the function
\eqref{mult-F-2} vanishes on the last points.

Since the discrete exponential functions $e_m(s)$ satisfy the
equality $e_m(s)=e_{m+N}(s)$, we do not need to consider them for all
values $m\in {\mathbb Z}^{\ge 0}$. It is enough to consider them for
$m\in \{ 1,2,\dots,N\}$.
By $\hat D_N^n$ we denote the set of integer $n$-tuples ${\bf m}=(m_1,m_2,\dots,m_n)$
such that
\[
 N\ge m_1>m_2>\cdots>m_n>0.
\]

We need a scalar product in the space of linear combinations of
the functions \eqref{mult-F-1}. It is natural to give it by the formula
\begin{equation}\label{mult-F-3a}
\langle e_{\bf m}({\bf s}),e_{{\bf m}'}({\bf s})   \rangle
\equiv \prod_{i=1}^n \langle e_{m_i}(s_i), e_{m'_i}(s_i)\rangle
:=\prod_{i=1}^n \sum_{s_i\in F_N}e_{m_i}(s_i)\overline{e_{m'_i}(s_i)}
=\delta_{{\bf m}{\bf m}'}
\end{equation}
where $m_i,m_i'\in \{ 1,2,\dots, N\}$. Here we used the relation \eqref{f-F-2}.
\medskip

\noindent
{\bf Proposition 1.} {\it For ${\bf m},{\bf m}'\in \hat D^n_N$
the discrete functions \eqref{mult-F-2}
 satisfy the orthogonality relation
\begin{equation}\label{mult-F-3}
\langle \tilde E^-_{\bf m}({\bf s}),
\tilde E^-_{{\bf m}'}({\bf s})\rangle =|S_n|
\sum_{{\bf s}\in \hat F_N^n} \tilde E^-_{\bf m}({\bf s})
\overline{\tilde E^-_{{\bf m}'}({\bf s})}
=\delta_{{\bf m}{\bf m}'},
 \end{equation}
where the scalar product is determined by formula \eqref{mult-F-3a}.}
\medskip

\noindent
{\bf Proof.} Since $m_1>m_2>\cdots>m_n>0$ and $m'_1>m'_2>\cdots>m'_n>0$, then due to
the definition of the scalar product we have
\[
\langle \tilde E^-_{\bf m}({\bf s}),
\tilde E^-_{{\bf m}'}({\bf s})\rangle =
\sum_{{\bf s}\in F_N^n} \tilde E^-_{\bf m}({\bf s})
\overline{\tilde E^-_{{\bf m}'}({\bf s})}
\]
\begin{equation}\label{mult-F-30}
=|S_n|^{-1}\sum_{w\in S_n} \prod_{i=1}^n
\sum_{s_i\in F_N} e_{m_{w(i)}}(s_i)\overline{e_{m'_{w(i)}}(s_i)}
= \delta_{{\bf m}{\bf m}'},
\end{equation}
where $(m_{w(1)},m_{w(2)},\dots ,m_{w(n)})$ is obtained from $(m_1,m_2,\dots,
m_n)$ by action by the permutation $w\in S_n$. Since functions
$\tilde E^-_{\bf m}({\bf s})$ are antisymmetric
with respect to $S_n$, then
\[
 \sum_{{\bf s}\in F_N^n} \tilde E^-_{\bf m}({\bf s})
\overline{\tilde E^-_{{\bf m}'}({\bf s})} =|S_n|
\sum_{{\bf s}\in \hat F_N^n} \tilde E^-_{\bf m}({\bf s})
\overline{\tilde E^-_{{\bf m}'}({\bf s})} .
\]
This proves the proposition.
 \medskip

Let $f$ be a function on $\hat F^n_N$ (or an antisymmetric function on $F_N^n$).
Then it can be expanded in the functions \eqref{mult-F-2} as
\begin{equation}\label{mult-F-4}
f({\bf s})=\sum_{{\bf m}\in \hat D_N^n}a_{\bf m}
\tilde E^-_{\bf m}({\bf s}).
 \end{equation}
The coefficients $a_{\bf m}$ are determined by the formula
\begin{equation}\label{mult-F-5}
a_{\bf m}=|S_n|\sum_{{\bf m}\in \hat F_N^n} f({\bf s})
\overline{\tilde E^-_{\bf m}({\bf s})}.
 \end{equation}
We have taken into account the facts that numbers of elements in
$\hat D_N^n$ and in $\hat F_N^n$
are the same and that the discrete functions \eqref{mult-F-2}
are orthogonal with respect to the scalar product
\eqref{mult-F-3}. We call expansions \eqref{mult-F-4} and \eqref{mult-F-5}
the {\it antisymmetric multivariate discrete Fourier transforms}.
These expansions can be written in terms of the exponential function
$E^-_{\bf m}({\bf s})=\det \left( \exp (2\pi {\rm i}m_is_j)\right)_{i,j=1}^n$,
\begin{equation}\label{mult-F-6}
f({\bf s})=N^{-n/2}\sum_{{\bf m}\in \hat D_N^n}a_{\bf m}
 E^-_{\bf m}({\bf s}),\ \ \
a_{\bf m}=N^{-n/2}|S_n|\sum_{{\bf m}\in \hat F_N^n} f({\bf s})
\overline{E^-_{\bf m}({\bf s})}.
 \end{equation}

\subsection{Symmetric multivariate discrete
Fourier transforms}
Let us give a symmetric multivariate discrete Fourier transforms.
For this we take the multivariate exponential functions \eqref{mult-F-1}
for integers $m_i$ such that
\[
N\ge m_1\ge m_2\ge \cdots\ge m_n\ge 1
\]
and make a symmetrization. We
obtain a finite version of the symmetric exponential functions
\eqref{det-SS},
\begin{equation}\label{mult-F-10}
\tilde E^+_{\bf m}({\bf s}):=
|S_n|^{-1/2}  {\det}^+ (e_{m_i}(s_j))_{i,j=1}^n=
|S_n|^{-1/2}N^{-1/2} E^+_{\bf m}({\bf s}),
 \end{equation}
where the discrete functions $e_m(s)$ are given by
\eqref{mult-F-1-1-1}.

The $n$-tuples ${\bf s}$ in \eqref{mult-F-10} run over $F_N^n\equiv
F_N\times \cdots\times F_N$ ($n$ times).
We denote by $\breve F_N^n$ the subset of $F_N^n$ consisting of ${\bf s}=(s_1,s_2,\dots,s_n)\in
F_N^n$ such that
\[
 s_1\ge s_2\ge \cdots \ge s_n.
\]
The set $\breve F_N^n$ is a finite subset of the closure of the fundamental
domain $F(S^{\rm aff}_n)$.

Note that acting by permutations $w\in S_n$ upon $\breve F_N^n$ we obtain
the whole set $F_N^n$, where each point, having some coordinates
$m_i$ coinciding, are repeated several times. Namely, a point ${\bf s}$ is contained
$|S_{\bf s}|$ times in $\{ w\breve F_N^n; w\in S_n\}$, where $S_{\bf s}$ is
the subgroup of $S_n$ consisting of elements $w\in S_n$ such that $w {\bf s}={\bf s}$.

By $\breve D_N^n$ we denote the set of integer $n$-tuples ${\bf m}=(m_1,m_2,\dots,m_n)$
such that
\[
 N\ge m_1\ge m_2\ge \cdots\ge m_n\ge 1.
\]

\noindent
{\bf Proposition 2.} {\it For ${\bf m},{\bf m}'\in \breve D_N^n$
the discrete functions \eqref{mult-F-10}
satisfy the orthogonality relation}
\begin{equation}\label{mult-F-11}
\langle \tilde E^+_{\bf m}({\bf s}),
\tilde E^+_{{\bf m}'}({\bf s})\rangle =
|S_n|\sum_{{\bf s}\in \breve F_M^n} |S_{\bf s}|^{-1}
\tilde E^+_{\bf m}({\bf s})
\overline{ \tilde E^+_{{\bf m}'}({\bf s})}
=|S_{\bf m}| \delta_{{\bf m}{\bf m}'}.
 \end{equation}

{\bf Proof.} This proposition is proved in the same way as Proposition 1, but
we have to take into account the difference between $\breve F_M^n$ and
$\hat F_M^n$.
Due to
the definition of the scalar product we have
\begin{alignat*}{2}
\langle \tilde E^+_{\bf m}({\bf s}),
\tilde E^+_{{\bf m}'}({\bf s})\rangle =&\;
\sum_{{\bf s}\in F_N^n} \tilde E^+_{\bf m}({\bf s})
\overline{\tilde E^+_{{\bf m}'}({\bf s})} \\
 =&\; |S_n|^{-1}|S_{\bf m}| \sum_{w\in S_n} \prod_{i=1}^n
\sum_{s_i\in F_N} e_{m_{w(i)}}(s_i)\overline{e_{m'_{w(i)}}(s_i)}
= |S_{\bf m}| \delta_{{\bf m}{\bf m}'},
\end{alignat*}
Here we have taken into
account that there appear additional summands (with
respect to \eqref{mult-F-30}) because some
summands on the right hand side of \eqref{mult-F-10} can coincide.

Since functions $\tilde E^+_{\bf m}({\bf s})$ are
symmetric with respect to $S_n$, then
\[
 \sum_{{\bf s}\in F_N^n} \tilde E^+_{\bf m}({\bf s})
\overline{\tilde E^+_{{\bf m}'}({\bf s})} =|S_n|
\sum_{{\bf s}\in \breve F_N^n} |S_{\bf s}|^{-1}
\tilde E^+_{\bf m}({\bf s})
\overline{\tilde E^+_{{\bf m}'}({\bf s})} ,
\]
where we have taken into account that under an action by $S_n$ upon
$\breve F_N^n$ a point ${\bf s}$ appears $|S_{\bf s}|$ times
in $F^n_N$.
This proves the proposition.
 \medskip

Let $f$ be a function on $\breve F^n_N$ (or a symmetric function on $F_N^n$).
Then it can be expanded in functions \eqref{mult-F-10} as
\begin{equation}\label{mult-F-12}
f({\bf s})=\sum_{{\bf m}\in \breve D_N^n}a_{\bf m}
\tilde E^+_{\bf m}({\bf s}).
 \end{equation}
The coefficients  $a_{\bf m}$ are determined by the formula
\begin{equation}\label{mult-F-13}
a_{\bf m}=|S_n| |S_{\bf m}|^{-1}\sum_{{\bf s}\in \breve F_N^n}
|S_{\bf s}|^{-1}f({\bf s}) \overline{\tilde E^+_{\bf m}({\bf s})}.
 \end{equation}

The expansions \eqref{mult-F-12} and \eqref{mult-F-13}
follows from the facts that numbers of elements in $\breve D_N^n$ and in $\breve F_N^n$
are the same and from the orthogonality relation
\eqref{mult-F-11}. We call expansions \eqref{mult-F-12} and \eqref{mult-F-13}
the {\it symmetric multivariate discrete Fourier transforms}.

\section{Eigenfunctions of (anti)symmetric Fourier transforms}
Let $H_n(x)$, $n=0,1,2,\dots$, be the well-known Hermite polynomials of one
variable. They satisfy the relation
 \begin{equation}\label{Herm-1}
 \int_{-\infty}^\infty e^{2\pi {\rm i}px}e^{-\pi p^2}
 H_m(\sqrt{2\pi}p)dp={\rm
i}^{-m}e^{-\pi x^2}H_m(\sqrt{2\pi}x)
 \end{equation}
(see, for example, subsection 12.2.4 in \cite{KVII}).

We create polynomials of many
variables
\begin{equation}\label{Her-m}
H_{\bf m}({\bf  x})\equiv H_{m_1,m_2,\dots,m_n}(x_1,x_2,\dots
,x_n):=H_{m_1}(x_1)H_{m_2}(x_2)\cdots H_{m_n}(x_n).
 \end{equation}
The functions
\begin{equation}\label{Her-2}
 e^{-|{\bf  x}|^2/2} H_{\bf m}({\bf  x}),\ \ \ \
m_i=0,1,2,\dots,\ \ \ \ i=1,2,\dots,n,
 \end{equation}
where $|{\bf  x}|$ is the length of the vector $x$, form an
orthogonal basis of the Hilbert space $L^2(\mathbb{R}^n)$
with the scalar product
$\langle f_1,f_2\rangle:=\int_{\mathbb{R}^n}
f_1(\mathbf{x})\overline{f_2(\mathbf{x})}d\mathbf{x}$,
where $d\mathbf{x}=dx_1\,dx_2\cdots dx_n$.

We make symmetrization and antisymmetrization of the functions
\[
{\mathcal H}_{\bf m}({\bf x}):=e^{-\pi |{\bf x}|^2} H_{\bf m}(\sqrt{2\pi}
{\bf x})
\]
(obtained from \eqref{Her-2} by replacing ${\bf x}$ by
$\sqrt{2\pi}{\bf x}$) by means of symmetric and
antisymmetric multivariate exponential functions:
\begin{equation}\label{Her-3}
\int_{\mathbb{R}^n} E^+_\lambda(\mathbf{x})
 e^{-\pi|{\bf  x}|^2} H_{\bf m}(\sqrt{2\pi}{\bf  x})=
 {\rm i}^{-|{\bf m}|} e^{-\pi|\lambda|^2}
 H^{\rm sym}_{\bf m}(\sqrt{2\pi}\lambda),
\end{equation}
\begin{equation}\label{Her-4}
\int_{\mathbb{R}^n} E^-_\lambda(\mathbf{x})
 e^{-\pi|{\bf  x}|^2} H_{\bf m}(\sqrt{2\pi}{\bf  x})
 ={\rm i}^{-|{\bf m}|} e^{-\pi|\lambda|^2}
 H^{\rm anti}_{\bf m}(\sqrt{2\pi}\lambda).
\end{equation}

It is easy to see that
the polynomials $H^{\rm sym}_{\bf m}$ and $H^{\rm anti}_{\bf m}$ indeed are
symmetric and antisymmetric, respectively, with respect to the
group $S_n$,
 \[
H^{\rm sym}_{\bf m}(w\lambda)=H^{\rm sym}_{\bf m}(\lambda),\ \ \ \
H^{\rm anti}_{\bf m}(w\lambda)=(\det w)H^{\rm anti}_{\bf m}
(\lambda),\ \ \ \ w\in S_n.
\]
For this reason, we may consider $H^{\rm sym}_{\bf m}(\lambda)$
for values of $\lambda=(\lambda_1,\lambda_2,\dots,\lambda_n)$
such that $\lambda_1\ge \lambda_2\ge \cdots \ge \lambda_n$ and
$H^{\rm anti}_{\bf m}(\lambda)$ for values of $\lambda$ such that
$\lambda_1>\lambda_2 >\cdots >\lambda_n$.
The polynomials $H^{\rm sym}_{\bf m}$ are of the form
 \begin{equation}\label{sym-H}
H^{\rm sym}_{\bf m}(\lambda)={\det}^+
\left(H_{m_i}(\lambda_j)\right)_{i,j=1}^n
 \end{equation}
and the polynomials $H^{\rm anti}_{\bf m}$ of the
form
 \begin{equation}\label{anti-H}
H^{\rm anti}_{\bf m}(\lambda) = \det
\left(H_{m_i}(\lambda_j)\right)_{i,j=1}^n.
 \end{equation}
Moreover, $H^{\rm anti}_{\bf m}(\lambda)=0$ if $m_i=m_{i+1}$ for some
$i=1,2,\dots,n-1$. For this reason, we may consider the polynomials
$H^{\rm sym}_{\bf m}(\lambda)$ for integer $n$-tuples ${\bf m}$ such that
$m_1\ge m_2\ge \cdots \ge m_n$ and the polynomials
$H^{\rm anti}_{\bf m}(\lambda)$ for integer $n$-tuples ${\bf m}$ such that
$m_1>m_2 >\cdots >m_n$.

Let us apply the symmetric Fourier transform \eqref{F-5} (we denote it as
$\mathfrak{F}$)
to the symmetric functions \eqref{sym-H}. Taking into
account formula \eqref{Her-3} we obtain
\begin{alignat*}{2}
\mathfrak{F}\left(  e^{-\pi|{\bf  x}|^2} H^{\rm sym}_{\bf m}
(\sqrt{2\pi}{\bf x})\right)=& \; \frac1{|S_n|} \int_{\mathbb{R}^n}
E^+_\lambda(\mathbf{x})
e^{-\pi|{\bf  x}|^2} H^{\rm sym}_{\bf m}
(\sqrt{2\pi}{\bf x}) d{\bf x}
 \notag\\
=& \; {\rm i}^{-|{\bf m}|} e^{-\pi|\lambda|^2}
 H^{\rm sym}_{\bf m}(\sqrt{2\pi}\lambda),
 \end{alignat*}
that is, {\it functions $e^{-\pi|{\bf  x}|^2} H^{\rm sym}_{\bf m}
(\sqrt{2\pi}{\bf x})$ are eigenfunctions of the symmetric Fourier
transform} $\mathfrak{F}$. Since these functions
for $m_i=0,1,2,\dots$, $i=1,2,\dots,n$, $m_1\ge m_2\ge \cdots \ge m_n$,
form an orthogonal basis of the Hilbert space $L_0^2(\mathbb{R}^n)$
of functions from $L^2(\mathbb{R}^n)$ symmetric with respect to $S_n$
(that is, of the Hilbert space $L^2(D_+)$),
then they constitute a complete set of
eigenfunctions of this transform. Thus, this transform has only four
eigenvalues ${\rm i}, -{\rm i}, 1, -1$ in $L_0^2(\mathbb{R}^n)$.
This means that we have $\mathfrak{F}^4=1$.

Now we apply the antisymmetric Fourier transform \eqref{F-3}
(we denote it as $\tilde{\mathfrak{F}}$) to the
antisymmetric function $e^{-\pi|{\bf  x}|^2} H^{\rm anti}_{\bf
m}(\sqrt{2\pi}{\bf x})$. Taking into account formula
\eqref{Her-4} we obtain
\begin{alignat*}{2}
\tilde{\mathfrak{F}}\left(  e^{-\pi|{\bf  x}|^2} H^{\rm anti}_{\bf
m}(\sqrt{2\pi}{\bf x})\right)=& \; \frac1{|S_n|} \int_{\mathbb{R}^n}
E^-_\lambda(\mathbf{x}) e^{-\pi|{\bf  x}|^2} H^{\rm anti}_{\bf
m}(\sqrt{2\pi}{\bf x}) d{\bf x}
 \notag\\
=& \; {\rm i}^{-|{\bf m}|} e^{-\pi|\lambda|^2}
 H^{\rm anti}_{\bf m}(\sqrt{2\pi}\lambda),
 \end{alignat*}
that is, {\it functions $e^{-\pi|{\bf  x}|^2} H^{\rm anti}_{\bf
m}(\sqrt{2\pi}{\bf x})$ are eigenfunctions of the
transform} $\tilde{\mathfrak{F}}$. Since these functions
for $m_i=0,1,2,\dots$; $i=1,2,\dots,n$,
$m_1>m_2>\cdots >m_n\ge 0$, form an
orthogonal basis of the Hilbert space $L_-^2(\mathbb{R}^n)$
of functions from $L^2(\mathbb{R}^n)$ antisymmetric with respect to
$W$, then they constitute a
complete set of eigenfunctions of this transform. Thus, this
transform has only four eigenvalues ${\rm i}, -{\rm i}, 1, -1$. This
means that, as in the previous case, we have
$\tilde{\mathfrak{F}}^4=1$.

\subsection*{Acknowledgements}

We are grateful for the hospitality extended to A.K. at the
Center de Recherches Math\'ematiques, Universit\'e de Montr\'eal,
during the preparation of this paper.
We acknowledge also partial support for this work from
the National Science and Engineering Research Council of
Canada, MITACS, the MIND Institute of Costa Mesa, California, and
Lockheed Martin, Canada.

\end{document}